\documentclass[11pt]{amsart}
\usepackage{amssymb}

\newcommand\Pf{{\noindent \bf Proof: }}
\newcommand\done{\hfill $\blacksquare$ \vspace{0.2 in}}
\newcommand\Z{\mathbb{Z}}
\newcommand\Q{\mathbb{Q}}
\newcommand\N{\mathbb{N}}
\def\P{{\mathcal{P}}}
\newcommand\calF{{\mathcal{F}}}
\newcommand\B{{\mathcal{B}}}
\newcommand\calO{{\mathcal{O}}}
\newcommand{\C}{{\mathfrak{C}}}
\newcommand{\D}{{\mathfrak{D}}}
\newcommand\Ftwo{\mathbb{F}_2}

\newcommand{\inflim}[1]{\displaystyle \lim_{#1 \rightarrow \infty}}
\newcommand{\invlim}{\displaystyle \lim_{\leftarrow}}
\newenvironment{Gal}{{\rm Gal\,}}{}
\newenvironment{Disc}{{\rm Disc \,}}{}
\newenvironment{Stab}{{\rm Stab \,}}{}
\newcommand{\abs}[1] {\left|#1\right|}
\newtheorem{theorem}{Theorem}[section]
\newtheorem{definition}[theorem]{Definition}
\newtheorem{lemma}[theorem]{Lemma}
\newtheorem{proposition}[theorem]{Proposition}
\newtheorem{corollary}[theorem]{Corollary}
\newtheorem{conjecture}[theorem]{Conjecture}

\title[density of prime divisors in arithmetic dynamics of quadratic polynomials]
{The density of prime divisors in the arithmetic dynamics of quadratic polynomials}
\author{Rafe Jones \\ University of Wisconsin-Madison}
\email{jones@math.wisc.edu}

\begin{document}


\begin{abstract}

Let $f \in \Z[x]$, and consider the recurrence given by $a_n = f(a_{n-1})$, with $a_0 \in \Z$.  Denote by $P(f,a_0)$ the set of prime divisors of this recurrence, i.e., the set of primes $p$ dividing some non-zero $a_n$, and denote the natural density of this set by $D(P(f,a_0))$.  The problem of determining $D(P(f,a_0))$ when $f$ is linear has attracted significant study, although it remains unresolved in full generality.   In this paper we consider the case of $f$ quadratic, where previously $D(P(f,a_0))$ was known only in a few cases.  We show $D(P(f,a_0)) = 0$ regardless of $a_0$ for four infinite families of $f$, including $f = x^2 +k$, $k \in \Z \backslash \{-1\}$.  The proof relies on tools from group theory and probability theory to formulate a sufficient condition for $D(P(f,a_0)) = 0$ in terms of arithmetic properties of the forward orbit of the critical point of $f$.  This provides an analogy to results in real and complex dynamics, where analytic properties of the forward orbit of the critical point have been shown to determine many global dynamical properties of a quadratic polynomial.  The article also includes apparently new work on the irreducibility of iterates of quadratic polynomials.

\end{abstract}

\maketitle

\section{Introduction}

Let $a_n = f(a_{n-1}, \ldots, a_{n-k})$ be a recursively defined sequence of order $k$, where $f$ is a polynomial in $k$ variables.  Many well-known sequences belong to this class, such as the Fibonacci numbers ($F_0 = F_1 = 1, F_{n} = F_{n-1} + F_{n-2}$), the Lucas numbers ($L_0 = 1, L_1 =  3, L_{n} = L_{n-1} + L_{n-2}$), the Fermat numbers $t_n = 2^{2^n}+1$ ($t_0=3, t_n = (t_{n-1} - 1)^2+1$), and the Mersenne numbers $m_n = 2^n - 1$ ($m_0 = 0, m_n = 2m_{n-1} + 1$).  The set of prime divisors of $a_n$, namely
\begin{equation} \label{setdef}
P(a_n) = \{ \text{$p$ prime : $p$ divides $a_i$ for some $i \geq 0$ with $a_i \neq 0$}\}
\end{equation}
has attracted great interest, although it often is too much to hope for a complete understanding of this set.  For instance, since the Fermat numbers $t_n$ are relatively prime, a full understanding of $P(t_n)$ would permit one to resolve the well-known conjecture that only finitely many of the $t_n$ are prime.  In this paper we consider the natural density $D(P(a_n))$ for certain $a_n$.  Recall that the natural density of a set $S$ of primes is
\begin{equation} \label{densedef}
D(S) = \inflim{x} \frac{\#\{p \in S : p \leq x\}}{\#\{p : p \leq x\}},
\end{equation}
provided that this limit exists.

When $a_n$ is a linear recurrence, the density of $P(a_n)$ has been much studied, e.g. \cite{ballot, hasse, lagarias, stephens}.  Note that $P(a_n)$ is infinite, excluding a few degenerate cases \cite{evertse, vanderpoorten}.  It is easy to see that $P(m_n)$ contains all $p \neq 2$ and $P(F_n)$ contains all $p$ (see e.g. \cite{vorobiev} for the latter).  However, the same is not true of $L_n$: Lagarias \cite{lagarias} shows $D(P(L_n)) = 2/3$.  Moreover, Hasse \cite{hasse} showed that the sequence $a_n = 2^n+1$, a relative of the Mersenne numbers, satisfies $D(P(a_n)) = 17/24$.  More generally, Stephens \cite{stephens} has studied the case $a_n = (b+1)a_{n-1} - ba_{n-2}.$  Assuming the Riemann Hypothesis for certain Kummer extensions, and with mild conditions on $a_0$ and $a_1$, he shows that
$D(P(a_n))$ is positive and depends in a delicate way on $b$.
Lagarias \cite{lagarias} remarks that under the Generalized Riemann Hypothesis it is possible to show that $D(P(a_n))$ exists and is positive for any non-degenerate second-order linear recurrence $a_n$.

We restrict our attention in this work to recurrences with first-order representations, that is, sequences $a_n = f(a_{n-1})$ with $f \in \Z[x]$ and $a_0 \in \Z$ (such a sequence also has a $k$th-order representation for any $k \geq 1$).   These sequences are orbits of integers in the arithmetic dynamical system defined by iteration of $f$.  We write $P(f, a_0)$ instead of $P(a_n)$, and our main concern is to determine $D(P(f, a_0))$ for various $f$ and $a_0$.  If $f$ is linear, then writing $f = bx + c$ we note that our sequence satisfies $a_n = (b+1)a_{n-1} - ba_{n-2},$ with $a_1 = ba_0 + c$.  This is the case treated in \cite{stephens}.

By contrast, the case where $\deg f \geq 2$ has received relatively little attention.  The main result is of Odoni \cite{odonigalit}, who proves that ``most" polynomials of degree $d$ have the property that all of their integer orbits contain a ``very small" proportion of the primes.  More specifically, he shows the following: let $A(d,N)$ denote the set of monic $f$ with $\deg f = d \geq 2$ and coefficients of absolute value at most $N$.  Then for any $\epsilon > 0$, as $N \rightarrow \infty$ the proportion of $f \in A(d,N)$ with $D(P(f, a_0)) < \epsilon$ for all $a_0 \in \Z$ approaches 1.
Odoni obtains his result by considering the generic monic of degree $d$ and applying a version of the Hilbert irreducibility theorem.  Thus the nature of the exceptional set is not explicit, and the result does not determine $D(P(f, a_0))$ for any specific $f$ and $a_0$.

On the other hand, in \cite{odoniwn} Odoni gives a thorough analysis of a sequence inspired by Euclid's proof of the infinitude of the primes, namely $w_0 = 2$, $w_n = 1 + w_0w_1 \cdots w_{n-1}$.  This sequence is also given by $w_0 = 2$, $w_n = f(w_{n-1})$, where $f = x^2 - x + 1$.  Odoni shows that
$D(P(f, a_0)) = 0$ for all $a_0 \in \Z$.

In this article we consider the case $\deg f = 2$.  We note that $P(f,a_0)$ is infinite as long as
$f \neq x^2$ and $\{f^n(a_0) : n \in \N\}$ is infinite; we prove this in Section \ref{inf}.  For $n \geq 1$ we denote the $n$th iterate of $f$ by $f^n$, and we take $f^0 = x$.
To give our primary results, we need a few preliminaries.  First, to deal with the case when some iterate of $f$ is reducible, it is desirable to have results for {\em translated} iterated sequences, that is, those of the form $g(f^n(a_0))$, where $g, f \in \Z[x]$.  We denote by $P(g,f, a_0)$ the set of primes dividing at least one non-zero term of this sequence.  Second, by a {\em rigid divisibility sequence}, \label{rdef} we mean a sequence $b_n$ of $S$-integers in $\Q$ with the property that for all $p \not\in S$, $v_p(b_n) > 0$ implies $v_p(b_{mn}) = v_p(b_n)$ for all $m \geq 1$.  This is a strengthening of the notion of divisibility sequence (i.e., $b_n \mid b_{mn}$ in the case $S = \emptyset$) and strong divisibility sequence (i.e., ${\rm gcd}\,(b_n, b_m) = b_{{\rm gcd}\,(n,m)}$ in the case $S = \emptyset$).  The cases of interest to us always have $S = \{2\}$ or $S = \emptyset$.

\begin{theorem} \label{mainone}
Let $f, g \in \Z[x]$ be monic with $f$ quadratic, and let $\gamma$ be the critical point of $f$.  Suppose that $g \circ f^n$ is irreducible for all $n \geq 0$ and the set $\{g(f^n(\gamma)\}$ is infinite.  If either
\begin{enumerate}
\item  $g$ divides $f^r$ for some $r \geq 0$ and the set $\{f^n(0) : n = 1, 2, \ldots\}$ is finite and does not contain $0$, or
\item the sequence $(g(f^n(\gamma)) : n = 1, 2, \ldots)$ is a rigid divisibility sequence,
\end{enumerate} then $D(P(g,f, a_0)) = 0$ for all $a_0 \in \Z$.
\end{theorem}

Our most frequent application of (1) is when $g(x) = x$, in which case $g$ divides $f^0 = x$.
The polynomials $f$ meeting the hypothesis of Theorem \ref{mainone} fall into several families,
and these are the subjects of our second main result.  We remark that the proof of Theorem \ref{mainone}, part 1 generalizes easily to the setting where $\Z$ is replaced by the ring of integers in a number field, and there it may apply to more polynomials (see p. \pageref{remkmainone}).

\begin{theorem} \label{maintwo}
Suppose one of the following holds:
\begin{enumerate}
\item $f = x^2 - kx + k$ for some $k \in \Z$
\item $f = x^2 + kx - 1$ for some $k \in \Z \backslash \{0,2\}$
\item $f = x^2+k$ for some $k \in \Z \backslash \{-1\}$
\item $f = x^2-2kx + k$ for some $k \in \Z \backslash \{\pm 1\}$.
\end{enumerate}
Then $D(P(f, a_0)) = 0$ for all $a_0 \in \Z$.
\end{theorem}

We note that the recurrence given by $w_n = k + w_0w_1 \cdots w_{n-1}$ satisfies
$w_n = f(w_{n-1})$, where $f = x^2 - kx + k$, whence part 1 of Theorem \ref{maintwo} generalizes the sequence considered in \cite{odoniwn}.  Note also that it is not true that $D(P(f,a_0)) = 0$ for all
monic, quadratic $f \in \Z[x]$.  Indeed, if $f = (x-k)^2 + k$ with $k = \pm p$ for any prime $p \geq 3$, then one can show $D(P(f,0)) = 1/3$.
However, it seems likely that $D(P(f,a_0)) = 0$ for nearly all $f$ whose critical point has infinite forward orbit (see Conjectures \ref{evstabconj} and \ref{bigconj}).  The ineffectivity of Theorems \ref{mainone} and \ref{maintwo} is due to the use of Siegel's theorem on finiteness of integral points on elliptic curves.

The method of proof of Theorems \ref{mainone} and \ref{maintwo} revolves around
a study of the Galois tower generated by the polynomials $g \circ f^n$ for $n = 1, 2, \ldots$.  We focus on the Galois groups $H_n(f,g)$ of the extensions $\Q(\text{roots of $g \circ f^n$}) / \Q(\text{roots of $g \circ f^{n-1}$})$.  To accomplish this, we amplify the techniques introduced in \cite{galmart} from the theory of stochastic processes.  This is the primary focus of sections \ref{evgalmart} and \ref{density}, which build up to a proof of Theorem \ref{mainone}.  One important step is Lemma \ref{zerolem}, which shows that to establish $D(P(f,a_0)) = 0$, it is generally enough to prove that $H_n(f,g)$ is as large as possible for infinitely many $n$.  By way of comparison, note that in \cite{odonigalit} and \cite{odoniwn} the results are achieved by showing $H_n(f,g)$ is as large as possible for all $n$.  To apply Theorem \ref{mainone} to specific families, we require results on the irreducibility of iterates of quadratic polynomials over $\Q$.  This is the subject of section \ref{stable}, and these results are apparently new.  In section \ref{families} we draw on the previous four sections to prove Theorem \ref{maintwo}.

A key step in the proof of Theorem \ref{mainone} is Theorem \ref{maxthm}, which shows that $H_n(f,g)$ is as large as possible only when the set $\{g(f^i(\gamma)) : i = 1, 2, \ldots, n\}$ satisfies certain arithmetic properties, where $\gamma$ is the critical point of $f$.  Thus the density of $P(f, a_0)$ depends on arithmetic properties of the critical orbit $\{g(f^i(\gamma)) : i = 1, 2, \ldots \}$ of $f$.
This makes for a striking analogy with complex and real dynamics, where analytic properties of the critical orbit of a quadratic polynomial have been shown to determine fundamental dynamical behavior of the polynomial.  For instance, if $f \in \mathbb{C}[z]$ is quadratic, then membership in the Mandelbrot set -- and equivalently the connectedness of the filled Julia set of $f$ -- is determined by whether the critical orbit remains bounded \cite[Section 3.8]{Devaney}.  If $f \in \mathbb{R}[x]$ is of the form $a - x^2$ for $-1/4 \leq a \leq 2$ the Collet-Eckmann condition states that if $\abs{f'}$ grows exponentially along the critical orbit, then $f$ exhibits stochastic behavior on an invariant interval.  An important question is whether almost all maps in this family satisfy the Collet-Eckmann condition; this question was recently resolved in the affirmative by Avila and Moreira \cite{avila}.

\section{Galois Processes and Eventual Martingales} \label{evgalmart}

Our main technique is to relate the density of $P(g,f, a_0)$ to properties of the Galois groups of $g \circ f^n$ over $\Q$.  In order to do this, we make use of the notion of {\em Galois processes} \cite{galmart}.

Let $K$ be a field and $g, f \in K[x]$ such that $g \circ f^n$ is separable for all $n \geq 0$ (we take $f^0$ to be the identity map).  Denote the splitting field of $g \circ f^n$ over $K$ by $K(g \circ f^n)$.  Clearly $K(g \circ f^n) \supseteq K(g \circ f^{n-1})$.  We now let $GP(f,g)$ be the Galois process of this tower, as defined in \cite[page~13]{galmart}; for completeness we give the construction in the present case.
Define $G_n(f,g)$ to be $\Gal(K(g \circ f^n)/K)$ and $H_n(f,g)$ to be $\Gal(K(g \circ f^n)/K(g \circ f^{n-1}))$.  Let $G(f,g) = \invlim G_n(f,g)$, and take $\P$ to be a Harr measure on $G(f,g)$ with $\P(G(f,g)) = 1$.  Letting $\B$ be the Borel sigma algebra, the triple
$(G(f,g), \P, \B)$ is then a probability space.   Denote by $\psi_n$ the natural projection $G(f,g)\rightarrow G_n(f,g)$, and define random variables $X_n$ on $G(f,g)$ as follows:
$$X_n(g) = \text{number of roots of $g \circ f^n$ fixed by $\psi_n(g)$}.$$
The data $(G(f,g), \P, \B, \{X_n\}_{n \geq 0})$ by definition
give a stochastic process, which we call the Galois process of $g \circ f^n$, and denote $GP(f,g)$.  Intuitively, this process resembles a
random walk through successively larger quotients of $G(f,g)$ that is compatible with the natural restriction maps.  Positions at each level are assigned a value based on the number of
roots of $g \circ f^n$ left fixed.  Note that it follows from basic properties
of Haar measure that
\begin{multline} \label{stronggoal}
\P(X_0 = t_0, \ldots, X_n = t_n) = \\ \frac{1}{\# G_n(f,g)}
\# \left\{\sigma \in G_n(f,g) : \mbox{$\sigma$ fixes $t_i$ roots of $g \circ f^i$ for
$i = 0,1, \ldots, n$} \right\}.
\end{multline}

The following result links $D(P(f, a_0))$ to Galois processes.  We wish to consider $f$ that may have reducible iterates, and thus it becomes necessary to deal with translates of recurrences.  Given a recurrence $a_n$ with terms in $ K$ and $g \in K[x]$, we refer to the sequence $g(a_n)$ as a {\em polynomial translate} of $a_n$.
\begin{theorem} \label{upbound}
Let $f,g \in \Z[x]$ be polynomials with $g \circ f^n$ separable for all $n$.  Let $a_n = f^n(a_0)$ with $a_0 \in \Z$
The density of prime divisors of $g(a_n) : n = 0,1,2,\ldots$ is bounded above by
$$\inflim{n} \P(X_n > 0),$$
where $X_n$ is the $n$th random variable in $GP(f,g)$.
\end{theorem}

{\bf \noindent  Remark}: It is also true that if $g = id$, then $\inflim{n} \P(X_n > 0)$ furnishes an upper bound for the density of primes $p$ such that $0$ is periodic in $\Z/ p\Z$ under iteration of $f$.  This follows from
the fact that $0$ is periodic in $\Z/ p\Z$ if and only if $f^{-n}(0) \cap \Z$ is non-empty for all $n \geq 1$; cf
\cite[Proposition 3.1]{galmart}.  Thus Theorem \ref{maintwo} shows that for the families in question, the density of $p$ with $0$ periodic in $\Z/ p\Z$ under $f$ is zero.

\smallskip

{\bf \noindent Remark}: With trivial modifications to the proof the Theorem can be made to apply with $\Z$ replaced by the ring of integers $\calO_K$ in a number field $K$.

\Pf Let
$$\Omega_N = \{p: p \nmid \Disc (g \circ f^N) \; \text{and $g(f^N(x)) \equiv 0 \bmod{p}$ has no solution in $\Z$}\}.$$
If $p \in \Omega_N$, then clearly $g(f^n(a_0)) \not\equiv 0 \pmod{p}$ for all $n \geq N$.  Thus
$p \nmid g(a_n)$ for all $n \geq N$.  There are only finitely many $p$ with
$p \mid g(a_n)$ for some $n < N$, and also only finitely many dividing $\Disc (g \circ f^N)$.  Therefore
\begin{equation} \label{krog}
D(\Omega_N) \leq D(\{p : \text{$p \nmid g(a_n)$ for all $n \geq 1$}\}).
\end{equation}

If $p \nmid \Disc g \circ f^N$, then $p$ cannot divide the discriminant of the splitting field $\Q(g \circ f^N)$
of $g \circ f^N$ \cite[Corollary 2, p. 157]{narkiewicz}, so $p$ is unramified in $\Q(g \circ f^N)$.  Now
$g(f^N(x)) \equiv 0 \pmod{p}$ having a solution in $\Z$ is equivalent to $g(f^N(x))$ having at least one linear factor in $(\Z / p\Z)[x]$.  Except for possibly finitely many $p$, this implies that
$p\calO_K = \mathfrak{P}_1 \cdots \mathfrak{P}_r$, where $K/\Q$ is obtained by adjoining a root of $g \circ f^N$, $\calO_K$ is the ring of integers, and at least one of the $\mathfrak{P}_i$ has residue class degree one \cite[Theorem 4.12]{narkiewicz}. This is equivalent to the disjoint cycle decomposition of the Frobenius conjugacy class at $p$ having a fixed point (in the natural permutation representation of $G_N(f,g) = \Gal (\Q(g \circ f^N)/\Q)$ on the roots of $g \circ f^N$).  From the Chebotarev Density Theorem it follows \cite[Proposition 7.15]{narkiewicz} that the density of $p$ with $p\calO_K$ having such a decomposition is
$$\frac{1}{\#G_N(f,g)} \#\{\sigma \in G_N(f,g) : \text{$\sigma$ fixes at least one root of $g \circ f^N$}\},$$
and by \eqref{stronggoal} this is simply $\P(X_N > 0)$.
Thus $\P(X_N > 0)  = 1 - D(\Omega_N)$.  From \eqref{krog}, we now have
$$\P(X_N > 0) \ \  \geq \ \ 1 - D(\{p : \text{$p \nmid g(a_n)$ for all $n \geq 1$}\}),$$
and this last expression is the same as $D(\{p : \text{$p \mid g(a_n)$ for some $n \geq 1$}\}).$ \done

In light of this result, we study the long-term behavior of $GP(f,g)$; this occupies the remainder of this section and also the next section.  Our main results show that $GP(f,g)$ has certain convergence properties, and require the following definition:
\begin{definition} \label{martdef}
A stochastic process $X_0, X_1, X_2, \ldots$ taking values in $\mathbb{Z}$ is a {\em martingale} if for all $n \geq 1$ and any $t_i \in \mathbb{Z}$,
$$E(X_n \mid X_{0} = t_{0}, X_1 = t_1, \ldots, X_{n-1} = t_{n-1}) = t_{n-1}.$$
We call $X_0, X_1, X_2, \ldots$ an {\em eventual martingale} if for some $N$ the process
$X_{N}, X_{N + 1}, X_{N + 2}, \ldots $ is a martingale.
\end{definition}

Martingales are important chiefly because they often converge in the following sense:
\begin{definition}
Let $X_1, X_2, \ldots$ be a stochastic process defined on the probability space
$(\Omega, \P, \calF)$.  The process {\em converges} if
$$\P \left(\omega \in \Omega : \text{$\inflim{n} X_n(\omega)$ exists} \right) = 1.$$
\end{definition}
Note that convergence does not depend on any finite collection of $X_n$.  We give one
martingale convergence theorem (see e.g. \cite[Section 12.3]{grimmett} for a proof)
\begin{theorem} \label{martconv}
Let $M = (X_1, X_2, \ldots)$ be a martingale whose random variables take nonnegative real values.  Then $M$ converges.
\end{theorem}
Since the random variables in $GP(f,g)$ take nonnegative integer values, we immediately have the following:
\begin{corollary} \label{evconstcor}
Suppose that $GP(f,g)$ is an eventual martingale.  Then
$$\P(\{\sigma \in G(f,g) : \text{$X_1(\sigma), X_2(\sigma), \ldots$ is eventually constant}\}) = 1.$$
\end{corollary}

A fruitful point of view of the groups $G_n(f,g)$ is as automorphism groups of rooted forests.  A rooted forest is a disjoint union of trees, each with a distinguished vertex called its root.  There is a natural filtration of the forest provided by the distance to a root.  We say a rooted forest is {\em complete of degree $d$} if each vertex is connected to precisely $d$ others, and that the forest has height $h$ if the length of the longest path is $h$.     We may form a rooted forest from the polynomials $g \circ f^n$ as follows:  take the roots (in the polynomial sense) of $g$ to be $V_0$, the set of roots (in the graph sense) of the forest, and for each $n \geq 1$, let $V_n$ be the set of roots of $g \circ f^n$.  Two vertices $v_{n-1} \in V_{n-1}$ and $v_n \in V_n$ are connected if $f(v_n) = v_{n-1}$.  For each root $\alpha$ of $g$, the set $\{v : \text{$f^i(v) = \alpha$ for some $i$}\}$ forms a tree with root $\alpha$.  It is easy to check that this defines a rooted forest.

A group $G$ acts on a rooted forest with vertices $V$ if there is an action of $G$ on the set $V$ such that
for all $\sigma \in G$, $v$ is connected to $v'$ if and only if $\sigma(v)$ is connected to $\sigma(v')$.

\begin{lemma} \label{stabav}
Let $G$ be a group acting on a rooted forest of height $n$, and denote the vertices of distance $i$ from a root by $V_i$.  Let $v_0 \in V_{n-1}$, and suppose the subgroup $H = \{\sigma \in G : \text{$\sigma(v) = v$ for all $v \in V_{n-1}$}\}$ acts transitively on  $E_{v_0} = \{v \in V_n : \text{$v$ is connected to $v_0$}\}$.  Finally, let $\sigma \in G$ satisfy
$\sigma(v_0) = v_0$.  Then
\begin{equation}
\frac{1}{\#H} \sum_{\tau \in H} \# \{v \in E_{v_0} : \sigma \tau(v) = v\} \\ = \\ 1.
\end{equation}
\end{lemma}

\Pf Since $\sigma (v_0) = v_0$, each element of $\sigma H$ permutes $E_{v_0}$.  Note that
$$\sum_{\tau \in H} \# \{v \in E_{v_0} : \sigma \tau(v) = v\} = \sum_{\tau \in H} \sum_{v \in E_{v_0}} \epsilon(v,\tau),$$
where $\epsilon(v,\tau) = 1$ if $\sigma \tau(v) = v$ and $\epsilon(v,\tau) = 0$ otherwise.  Interchanging the sums yields
\begin{equation} \label{ptstab}
\sum_{v \in E_{v_0}} \# \Stab_{\sigma H} (v).
\end{equation}
Now given $v \in E_{v_0}$, by the transitivity of $H$ we can choose $\tau' \in H$ with $\tau' (v) = \sigma^{-1}(v)$.  We then have $\sigma \tau (v) = v$ if and only if $\tau'^{-1}\tau(v) = v$.  Thus
$$\# \Stab_{\sigma H} (v) = \#\{\tau \in H : \tau'^{-1}\tau(v) = v\} = \# \Stab_H (v).$$
By standard group theory and the transitivity of $H$, we have $\# \Stab_H (v) = \#H / \# E_{v_0}$.  Hence the expression in \eqref{ptstab} is $\#H$.  \done

\begin{theorem}  \label{evmartthm}
Let $f, g \in K[x]$ be such that $g \circ f^n$ is separable for all $n \geq 0$.  Suppose that there exists $n_0$ such that for all $n \geq n_0$ and every root $\alpha$ of $g \circ f^{n-1}$, the polynomial $f(x) - \alpha$ is irreducible over the splitting field $K(g \circ f^{n-1})$ of $g \circ f^{n-1}$.  Then $GP(f, g)$ is an eventual martingale.
\end{theorem}

\Pf By the definition of eventual martingale and conditional expectation, we must show that
there is an $N$ such that for $n > N$, the expression
\begin{equation} \label{sloop1}
\sum_k k \cdot \frac{\P (X_{N} = t_N, \ldots, X_{n-1} = t_{n-1}, X_n = k )}
{\P (X_N = t_N, \ldots, X_{n-1} = t_{n-1})}
\end{equation}
is equal to $t_{n-1}$.  We take $N$ to be $n_0$.
Put
\begin{eqnarray*}
S & = & \{\sigma \in G_n(f,g) : \text{$\sigma$ fixes $t_i$ roots of $g \circ f^i$ for $N
\leq i \leq n-1$}\} \\
S_k & = & \{\sigma \in S : \text{$\sigma$ fixes $k$ roots of $g \circ f^n$}\}
\end{eqnarray*}
From the basic property of GP$(f,g)$ given in \eqref{stronggoal},
the expression in \eqref{sloop1} is equal to
\begin{equation*} \label{sloop2}
\sum_k k \cdot \frac{\# S_k }
{\# S }
\end{equation*}
This in turn may be rewritten as
\begin{equation} \label{sloop4}
\frac{1}{\#S} \sum_{\sigma \in S} \left( \text{number of roots of
$g \circ f^n$ fixed by $\sigma$} \right).
\end{equation}

We denote $H_n(f,g)$ by $H$.  Note that if $\tau \in H$, then $\tau$ fixes the roots of
$g \circ f^i$ for $N \leq i \leq n-1$.  Thus $S$ is
invariant under multiplication by $H$, whence $S$ is a union of
cosets of $H$.   Let $\sigma H \subseteq S$, and let $T$ be complete binary infinite rooted forest corresponding to the translated iterates $g \circ f^n$.
Note that if $T_n$ denotes the height $n$ rooted forest consisting of the first $n$ levels
of $T$, then $G_n(f,g)$ acts as automorphisms of $T_n$ and $H$ is the subgroup stabilizing
the vertices of height $n-1$.

Let $\alpha_0$ be a vertex of height $n-1$, and suppose that $\sigma(\alpha_0) = \alpha_0$.
We have the following factorization of $g \circ f^n$ over $K(g \circ f^{n-1})$:
\begin{equation} \label{prood}
g \circ f^n = \prod_{\text{$\alpha$ root of $g \circ f^{n-1}$}} f(x) - \alpha.
\end{equation}
Thus the irreducibility of $f(x)-\alpha_0$ implies that $H$ acts transitively on the set
$E_{\alpha_0} = \{\text{$\beta$ root of $g \circ f^n$} : f(\beta) = \alpha_0\}$, which is the set of vertices of height $n$ connected to $\alpha_0$.   By Lemma \ref{stabav}, we have
$$
\frac{1}{\#H} \sum_{\tau \in H} \# \{\beta \in E_{\alpha_0} : \sigma \tau(\beta) = \beta\} \\ = \\ 1.
$$

On the other hand, if $\alpha$ is a vertex of height $n-1$ and $\sigma(\alpha) = \alpha'$ for some
$\alpha' \neq \alpha$ then $\sigma$ maps the roots of $f(x) - \alpha$ to the roots of
$f(x) - \alpha'$, and hence fixes no roots of $f(x) - \alpha$.  Since $\sigma \in S$, we have that $\sigma$ fixes $t_{n-1}$
roots of $f^{n-1}$, so it follows that
$$\frac{1}{\#H}\sum_{\tau \in H} \# \{\text{roots of $g \circ f^n$ fixed by $\sigma \tau$}\} = t_{n-1}.$$
Since $S$ is a union of cosets of $H$, the expression in
\eqref{sloop4} equals $t_{n-1}$.  \done

We remark that it follows from the factorization in \eqref{prood} that the extension $\Q(g \circ f^n)/\Q(g \circ f^{n-1})$ is the compositum of at most $\deg g \circ f^{n-1}$ quadratic extensions.
Thus $H_n(f,g) \cong (\Z / 2\Z)^m$ for some $0 \leq m \leq \deg g \circ f^{n-1}$.  If $m = \deg g \circ f^{n-1}$, then we call $H_n(f, g)$ {\em maximal}.  \label{maximal}

\begin{lemma} \label{disclem}
Let $f,g \in K[x]$ where $f$ has degree $d$ and leading coefficient $a$, and let $\gamma_1, \gamma_2, \ldots, \gamma_{d-1}$ be the critical points of $f$.
Put $\Delta_n = \Disc g \circ f^n$.  Then for all $n \geq 1$ we have
$$\Delta_n = \Delta_{n-1}^d((-1)^{d-1}da)^{d^n \cdot \deg g} \prod_{i=1}^{d-1} g(f^n(\gamma_i)).$$
\end{lemma}

\Pf Recall that the resultant $R(h_1, h_2)$ of two polynomials $h_1$ and $h_2$ is defined to be
$\prod_{h_1(\alpha) = 0} h_2(\alpha)$.
By definition, $\Delta_n = R(g \circ f^n, (g \circ f^n)')$.  Writing
$g \circ f^n$ as $(g \circ f^{n-1}) \circ f$, we have
$(g \circ f^n)' = f' \cdot (g \circ f^{n-1})' (f).$
Therefore $$R(g \circ f^n, (g \circ f^n)') = R(g \circ f^n, f') \cdot (R(g \circ f^{n-1}, (g \circ f^{n-1})'))^{\deg f},
$$ which gives $\Delta_n = R(g \circ f^n, f') \Delta_{n-1}^d$.  Now
\begin{multline*}
\prod_{g(f^n(\alpha)) = 0} f'(\alpha)  \; = \; \prod_{g(f^n(\alpha)) = 0} ad \prod_{i=1}^{d-1} (\alpha - \gamma_i) \; = \; \\ ((-1)^{d-1}ad)^{d^n \cdot \deg g}  \prod_{i=1}^{d-1} \prod_{g(f^n(\alpha)) = 0}  (\gamma_i - \alpha) \; = \; ((-1)^{d-1}ad)^{d^n \cdot \deg g}  \prod_{i=1}^{d-1} g(f^n(\gamma_i)).
\end{multline*}
This proves the lemma.  \done

One consequence of Lemma \ref{disclem} is that $g \circ f^n$ is separable for all $n$ if and only if $g$ is separable and $0$ is not contained in the forward orbit of any critical point of $f$.  We also remark that in the case $\deg f = 2$, Lemma \ref{disclem} gives
\begin{equation}
\Delta_n = \Delta_{n-1}^2(-2a)^{2^n \cdot \deg g} g(f^n(\gamma)),
\end{equation}
where $\gamma$ is the unique critical point of $f$.

We call a polynomial $f \in K[x]$ {\em critically finite} if all of its critical points have finite forward orbits, or equivalently are eventually periodic.  In particular, a quadratic polynomial $f = ax^2 + bx + c$ is critically finite if and only if the forward orbit of $-b/2a$ is finite.  We remark that it follows from Lemma \ref{disclem} and the appendix to this article that the compositum of $\Q(g \circ f^n), n = 1, 2, \ldots$ is finitely ramified if and only if $f$ is critically finite (cf \cite{hajir}).

We now specialize to the case where $f$ is quadratic and $K$ is a number field.

\begin{theorem} \label{monk}
Let $K$ be a number field, $f, g \in K[x]$ with $f$ quadratic, and suppose $G_0(f,g) = \Gal (K(g)/K)$ is a 2-group.  Suppose further that $f$ is critically infinite, and $g \circ f^n$ is irreducible for all $n \geq 0$.  Then $GP(f,g)$ is an eventual martingale.
\end{theorem}

{\bf \noindent Remark:} For the conclusion of Theorem \ref{monk} to hold, we only need the sequence
$\{g(f^n(\gamma)) : n = 1, 2, \ldots\}$ to contain only finitely many squares.  This happens for many critically finite polynomials, e.g. $x^2 - 2$ and $(x-k)^2 + k$ for all nonsquare $k \in \Z$.  \smallskip

\Pf  Let $K_n$ be the splitting field of $g \circ f^n$ over $K$,
Let $V_n$ be the roots of $g \circ f^n$, and let $\C$ be the partition of $V_n$ into the sets $\{\text{roots of $f(x) - \alpha$}\},$ where $\alpha$ varies over roots of $g \circ f^{n-1}$.  Let $\sigma_{\C} \in {\rm Sym} (V_n)$ be the permutation associated to
$\C$, i.e. the unique permutation whose orbits are precisely the sets belonging to $\C$.  We wish to show that $\sigma_{\C} \in G_n(f,g)$ for $n$ sufficiently large.

For all $n \geq 1$, it follows from \eqref{prood} that the field $K_n$ is obtained from $K_{n-1}$ by adjoining square roots.  Since $G_0(f,g)$ is a 2-group, it follows that $G_n(f,g)$ is a 2-group as well, and thus has non-trivial center.  Let $\delta$ be a nontrivial element of $Z(G_n(f,g))$, and $\D$ the corresponding central fiber system (see \cite[Proposition-Definition 4.10]{galmart}).  Then by \cite[Theorem 4.9]{galmart} we have
either $\sigma_{\C} = \delta$ or $G_n(f,g)$ is composed entirely of alternating permutations.

In the latter case, $\Delta_n = \Disc g \circ f^n$ is a square in $K$, and by Proposition \ref{disclem} this implies that $g(f^n(\gamma))$ is a square in $K$, where $\gamma$ is the critical point of $f$.  Now for $n \geq 2$, this last statement implies the curve $C : y^2 = g(f^2(x))$ has a pair
of points that are $S$-integers, indeed a pair with $x = f^{n-2}(\gamma)$ (we may take $S$ to be the primes in $\calO_K$ dividing either $2$ or the denominator of one of the coefficients of $f$).  Since $f$ is critically infinite, $f^k(\gamma) \neq f^j(\gamma)$ for all $j,k$, so that the pairs obtained from different $n$ are all distinct.  Since $\deg g \circ f^2 \geq 4$, and $g \circ f^2$ is irreducible and thus separable, the genus of $C$ is at least 1, whence by Siegel's Theorem \cite[p.353]{jhsdioph} $C$ has only finitely many $S$-integral points.  Thus $g(f^n(\gamma))$ is a square for only finitely many $n$.  Denote the largest such $n$ by $n'$, and take $n_0 = n'+1$.  Thus for $n \geq n_0$ we have shown that $\sigma_{\C} = \delta$, and it follows that
$\sigma_{\C} \in G_n(f,g)$.  Indeed, $\sigma_{\C} \in H_n(f,g) \stackrel{\rm def}{=} \Gal (K_n/K_{n-1})$, whence $H_n(f,g)$ is nontrivial.

Therefore by \eqref{prood} there is some root $\alpha$ of $g \circ f^{n-1}$ such that $f(x) - \alpha$ is irreducible over $K_{n-1}$.  Since $g \circ f^{n-1}$ is irreducible over $\Q$, the action of $G_{n-1}(f,g)$ is transitive on the roots of $g \circ f^{n-1}$, and it follows that $f(x) - \alpha$ is irreducible over $K_{n-1}$ for all
$\alpha$.  The irreducibility of $g \circ f^n$ over $K$ implies separability, and we now have from Theorem \ref{evmartthm} that $GP(f,g)$ is an eventual martingale. \done

\section{Density results} \label{density}

In this section we prove Theorem \ref{mainone}.  Note that by Theorem \ref{monk}, the hypotheses of Theorem \ref{mainone} imply that $GP(f,g)$ is an eventual martingale.  In order to parlay this into a zero-density result, we use the following:

\begin{lemma} \label{zerolem}
If $GP(f,g)$ is an eventual martingale and $H_n(f,g)$ is maximal for
infinitely many $n$, then $GP(f,g)$ converges to $0$, i.e.
$$\inflim{n} \P(X_n = 0) = 1.$$
\end{lemma}

\Pf  The proof of \cite[Lemma 5.2]{galmart}
requires only straightforward adaptation to apply to $GP(f,g)$.
The proof of \cite[Theorem 1.3]{galmart}
only uses \cite[Lemma 5.2]{galmart} and that the stochastic process is eventually constant with probability 1, which follows from Corollary \ref{evconstcor}.  \done

We now give two results that relate the maximality of $H_n(f,g)$ to the translated critical orbit
$\{g(f^n(\gamma)) : n \geq 1\}$.

\begin{lemma} \label{stoll}
Let $K$ be a number field and $f,g \in K[x]$ with $f = ax^2 + bx + c$, and let $\gamma$ be the critical point of $f$.  Suppose that $g \circ f^n$ is irreducible over $K$ for all $n \geq 1$.  Then for $n \geq 2$, $H_n(f,g)$ is maximal if and only if $g(f^n(\gamma))$ is not a square in $K(g \circ f^{n-1})$.
\end{lemma}

\Pf We require only minor modifications of \cite[Lemma 1.6]{stoll}.  Let $d = \deg g \circ f^{n-1}$,
$K_n = K(g \circ f^{n})$,  $K_{n-1} = K(g \circ f^{n-1})$, and denote the roots of $g \circ f^{n-1}$ by
$\alpha_j, j = 1, \ldots, d$.  From \eqref{prood}, we have that $K_n/K_{n-1}$ is a 2-Kummer extension obtained by adjoining to $K_{n-1}$ the square roots of $b^2 - 4ac + 4a\alpha_j$ for $j = 1, \ldots, d$.  Thus $[K_n : K_{n-1}] = 2^{d - \dim V}$, where  $V$ is the $\Ftwo$-vector space of all $\epsilon_1, \ldots, \epsilon_d \in \Ftwo^d$ such that $\prod_{j=1}^d  (b^2 - 4ac + 4a\alpha_j)^{\epsilon_j}$ is a square in
$K_{n-1}$.  By the irreducibility of $g \circ f^n$, we have that $\Gal (K_n/K)$ acts
transitively on the $\alpha_j$.  We now follow precisely the reasoning in the proof of \cite[Lemma 1.6]{stoll} to show that $V \neq 0$ if and only if $(1, \ldots, 1) \in V$.  Therefore $H_n(f,g)$ is maximal if and only if $\prod_{j=1}^{d}  b^2 - 4ac + 4a\alpha_j$ is a square in $K_{n-1}$.  Since $n \geq 2$, we have $d$ even, and thus
$H_n(f,g)$ is maximal if and only if
$$\prod_{j=1}^{d}  \alpha_j - \frac{b^2}{4} + c$$ is a square in $K_{n-1}$.  As in the proof of Proposition \ref{stabcrit}, this is equivalent to $g(f^n (\gamma))$ being a square in $K_{n-1}$.
\done

\begin{theorem} \label{maxthm}
Let $K$ be a number field with ring of integers $\calO_K$, and let $f,g \in \calO_K[x]$ with $f = ax^2 + bx + c$, and let $\gamma$ be the critical point of $f$.  Suppose that $g \circ f^n$ is irreducible over $K$ for all $n \geq 1$, and denote by $v_{\mathfrak{p}}$ the $\mathfrak{p}$-adic valuation for $\mathfrak{p}$ a prime in $\calO_K$.  If $n \geq 2$ and there exists $\mathfrak{p}$ with $v_{\mathfrak{p}}(g(f^n(\gamma)))$ odd, $v_{\mathfrak{p}}(g(f^m(\gamma))) = 0$ for all $1 \leq m \leq n-1$, and $v_{\mathfrak{p}}(2a) = 0$ then $H_n(f,g)$ is maximal.
\end{theorem}

\Pf  It follows from Lemma \ref{disclem} and \cite[Corollary 2, p. 157]{narkiewicz} that the only primes ramifying in $K(g \circ f^{n-1})$ are a subset of those dividing $2a$ and those dividing the numerator of $g(f^m(\gamma))$ for some $0 \leq m \leq n-1$.  Thus if
$v_{\mathfrak{p}}(g(f^m(\gamma))) = 0$ for all $1 \leq m \leq n-1$, and $v_{\mathfrak{p}}(2a) = 0$ then
$\mathfrak{p}$ is unramified in $K(g \circ f^{n-1})$.  Hence if $v_{\mathfrak{p}}(g(f^n(\gamma)))$ is odd, then $g(f^n(\gamma))$ cannot be a square in $K(g \circ f^{n-1})$.  By Lemma \ref{stoll} we then have
$H_n(f,g)$ maximal.  \done

As an illustration of Theorem \ref{maxthm}, consider the case of $K = \Q$, $f(x) = x^2+5$, $g(x) = x$.  We have $f(\gamma) = 5$,
$f^2(\gamma) = 2 \cdot 3 \cdot 5$, $f^3(\gamma) = 5 \cdot 181$, and
$f^4(\gamma) = 2 \cdot 3 \cdot 5 \cdot 23 \cdot 1187$.  In each case there is a $p \neq 2$ with
$v_p( f^n(\gamma)) = 1$ but $v_p(f^m(\gamma)) = 0$ for all $1 \leq m \leq n-1$, and therefore
$H_n(f)$ is maximal for $1 \leq n \leq 4$.

Theorem \ref{maxthm} says that the maximality of $H_n(f,g)$ can be ascertained by examining arithmetic properties of the critical orbit.  This makes for a striking analogy with complex and real dynamics, where analytic properties of the critical orbit of a quadratic polynomial have been shown to determine fundamental dynamical behavior of the polynomial.  See the last paragraph of the Introduction for more details.

We make one final remark before proving Theorem \ref{mainone}.  Consider $f \in \Z[x]$ quadratic, and let $g(x) = x$.  To analyze the arithmetic of the critical orbit of $f$, we split the primes dividing some
$f^n(\gamma)$ into two classes.  We say that $p$ is {\em isolated} if $v_p (f^n(\gamma)) > 0$ for some $n$ and there is no other $n'$ with $v_p (g(f^{n'}(\gamma))) > 0$.  When we reduce the coefficients of $f$ modulo $p$ and iterate, this behavior corresponds to $0$ being in the orbit of $\gamma$ in $\Z / p\Z$ but not in the periodic part.  If $p$ is not isolated, then there are infinitely many $n'$ with
$v_p (f^{n'}(\gamma)) > 0$, and we call $p$ {\em recurrent}.  This corresponds to $0$ being in the periodic part of the orbit of $\gamma$ in $\Z / p\Z$.   Heuristics and computation suggest that most primes $p$ dividing some element of the critical orbit are isolated, but in general there appear to be infinitely many recurrent primes, and their behavior seems difficult to control.  There is one exception, namely when all primes dividing an element of the translated critical orbit are recurrent.  This is the case when $\gamma = 0$, e.g. for $f(x) = x^2 + k$, $k \in \Z$.  At the other extreme are $f$ such that all but finitely many primes appearing in the critical orbit are isolated.  This is the case when
$\{f^n(0) : n \geq 1\}$ is finite and does not contain $0$, as in the case of $f = x^2 - kx + k$, for which $0$ maps to the fixed point $k$.  These two extremes correspond to parts (2) and (1), respectively, of Theorem \ref{mainone}.

\smallskip

{\bf Proof of Theorem \ref{mainone}}:
We apply Theorem \ref{monk} to get that $GP(f)$ is an eventual martingale.  Then by Lemma \ref{zerolem} and Theorem \ref{upbound} we need only show that $H_n(f,g)$ is maximal for infinitely many $n$.
Suppose first that we are in case 1, and let $\{f^n(0) : n = 1, 2, \ldots\} = \{b_1, \ldots, b_k\}$ and
$\{g(f^n(0)) : n = 1, 2, \ldots\} = \{\beta_1, \ldots, \beta_k\}$.  We remark that $g(f^n(0)) \mid f^{n+r}(0)$, so that each $\beta_i$ divides some $b_j$.
The outline of the argument is this:
it follows from Siegel's theorem on $S$-integral points that for infinitely many $n$, the rational number $\abs{g(f^n(\gamma))/d}$ is not a square for any $d \mid \prod_i b_i$.  We use this fact to show that
there must exist $p$ with $v_p(g(f^n(\gamma)))$ odd and $v_p(g(f^j(\gamma))) = 0$ for all $j < n$.  It follows from Theorem \ref{maxthm} that $H_n$ is maximal, and this proves the theorem.

To apply Siegel's theorem, we assume there are infinitely many $n$ with $g(f^n(\gamma)) > 0$; the argument is similar if there are infinitely many $n$ with $g(f^n(\gamma)) < 0$.  If $g(f^n(\gamma))/\abs{d}$ is a square for $n \geq 2$, then the curve $C_d: \abs{d}y^2 = g(f^2(x))$ has a pair of $S$-integral points ($S = \{2\}$)  with $x = f^{n-2}(\gamma)$.  By assumption the set $\{g(f^n(\gamma))\}$ is infinite, so all points obtained in this way are distinct.  Since $g(f^2(x))$ is irreducible and hence separable, it follows that $C_d$ has genus at least $1$.  From Siegel's Theorem on $S$-integral points \cite[p. 353]{jhsdioph} and the fact that there are only finitely many values for $d$, we have that for all $n$ large enough, $g(f^n(\gamma))/\abs{d}$ is not a square for any $d$.  Since $g(f^n(\gamma)) > 0$ for infinitely many $n$, we have for these $n$ that $\abs{g(f^n(\gamma))/d}$ is not a square for any $d$, as desired.

Suppose for a moment that $\gamma \in \Z$.  We now show that for all $n$ large enough, there exists
$d \mid \prod_i b_i$ such that $g(f^n(\gamma))/d$ is an integer not divisible by any primes
$p \mid \prod_i b_i$.  This is equivalent to showing that $v_p(g(f^n(\gamma))) \leq v_p(\prod_i b_i)$ for all $p \mid \prod_i b_i$.  Let $m = \max_i{v_p(b_i)}$, and
suppose $v_p(g(f^{n_0}(\gamma))) > m$.  By assumption $g \mid f^r$, say $gh = f^r$.  Then multiplying through by $h(f^{n_0}(\gamma))$ gives $f^{n_0+r}(\gamma) \equiv 0 \bmod{p^{m+1}}$.   Thus for $j > n_0 + r$ we have $f^j(\gamma) \equiv  b_i \not\equiv 0 \bmod{p^{m+1}}$
for some $i$.  Hence for these $j$ we have $g(f^j(\gamma)) \equiv  \beta_i \not\equiv 0 \bmod{p^{m+1}}$.
Therefore by taking $n$ large enough we may ensure $v_p(g(f^n(\gamma))) \leq m \leq v_p(\prod_i b_i)$ for all $p \mid \prod_i b_i$.   If $\gamma \not\in \Z$, then we wish to show that for all $n$ large enough there exists an {\em odd} $d \mid \prod_i b_i$ such that the numerator of $g(f^n(\gamma))/d$ is an integer not divisible by any primes $p \mid \prod_i b_i$ (the denominator of $g(f^n(\gamma))/d$ is an even power of $2$).  We proceed exactly as before, noting that the numerator of $g(f^n(\gamma))$ is odd, and thus $d$ can be taken to be odd as well.

From the previous two paragraphs it follows that for infinitely many $n$, there exists a prime
$p$, not dividing any of the $b_i$, such that $v_p(g(f^n(\gamma)))$ is odd.  Therefore writing $gh = f^r$ and multiplying by $h(f^n(\gamma))$ we have $f^{n+r}(\gamma) \equiv 0 \bmod p$, so for each $j > n + r$ we have $f^j(\gamma) \equiv b_i \not\equiv 0 \bmod{p}$ for some $i$.  If $v_p(g(f^l(\gamma))) > 0$ for some $l > n$, then multiplying by $h(f^l(\gamma))$ gives $f^{l+r}(\gamma) \equiv 0 \bmod p$, a contradiction.  It follows that $v_p(g(f^l(\gamma))) = 0$ for all $l < n$ as well, since otherwise $v_p(g(f^n(\gamma))) = 0$, a contradiction (hence $p$ is isolated in the terminology of the discussion preceding this proof).  It remains to show we can take $n$ large enough to ensure $p \neq 2$.  If $\gamma \not\in \Z$, then the numerator of $f^n(\gamma)$ is odd and for all $n$, whence the same is true for $g(f^n(\gamma))$, since $g$ divides $f^r$.  In this situation $v_2(g(f^n(\gamma))$ is never odd.  If $\gamma \in \Z$ then either $2 \mid \prod_i b_i$ or there is at most one $n$ with $v_2(g(f^n(\gamma))) > 0$.  Hence $p \neq 2$ for $n$ large enough, and  the theorem now follows from Theorem \ref{maxthm}.

Now suppose that we are in case 2, and let $c_n = g(f^n(\gamma))$ for $n \geq 1$.  Note that $(c_n)$ is a rigid divisibility sequence (see p. \pageref{rdef}), with $S = \{2\}$ if $\gamma \not\in \Z$ and $S = \emptyset$ if $\gamma \in \Z$.  Consider $c_q$, where $q$ is a prime.  If $v_p(c_q) > 0$, then
we have either $v_p(c_1) > 0$ or $v_p(c_n) = 0$ for $1 \leq n < q$.  If $v_p(c_1) > 0$, then $v_p(c_1) = v_p(c_q)$, whence $v_p(c_q/c_1) = 0$.

Thus we need only show that for infinitely many $q$ there exists $p \neq 2$ with $v_p(c_q/c_1)$ odd: by Theorem \ref{maxthm} this implies that $H_q(f,g)$ maximal for infinitely many $q$, which completes the proof.

Note that $v_2(c_q/c_1)$ cannot be odd for $q > 2$.  Indeed, if $\gamma \not\in \Z$ then
since $g$ and $f$ are monic, $v_2(c_n) = -2^n$ for all $n$.  If $\gamma \in \Z$ then either $v_2(c_1) > 0$, implying $v_2(c_n/c_1) = 0$ for all $n$, or $v_2(c_1) = 0$ and  $v_2(c_2) > 0$, implying $v_2(c_n) = 0$ for $n$ odd, or $v_2(c_1) = 0$ and $v_2(c_2) = 0$, implying that $v_2(c_n) = 0$ for all $n$.  The first two implications follow from the fact that $c_n$ is a rigid divisibility sequence, and the last since
$g(f(\gamma)) \equiv 1 \bmod{2}$ and $g(f^2(\gamma)) \equiv 1 \bmod{2}$ imply either $g(0) \equiv g(1) \equiv 1 \bmod{2}$ or $f^n(\gamma) \equiv f^m(\gamma) \bmod{2}$ for all $m, n$.

Let $q > 2$.  To show there exists $p \neq 2$ with $v_p(c_q/c_1)$ odd, it is thus enough to show that $\abs{c_q/c_1}$ is not a square in $\Q$.  Suppose that infinitely many of the $c_q$ are positive; the proof is similar in the case that infinitely many are negative.   If $q \geq 2$ and $c_q/\abs{c_1}$ is a square, then the curve $C: \abs{t_1}y^2 = f^2(x)$ has a pair of integral points with $x = f^{q-2}(\gamma) = c_{q-2}$.  Since $\{c_n : n \geq 1\}$ is infinite by assumption, each pair obtained in this way is distinct.  Since $f^2$ is separable, it follows that $C$ has genus 1, whence by Siegel's Theorem \cite[p. 353]{jhsdioph} there can be only finitely many such pairs, and thus we have $c_q/\abs{c_1}$ not a square for $q$ sufficiently large.  By assumption $c_q = \abs{c_q}$ for infinitely many $q$, and this completes the proof.  \done

{\bf \noindent Remark:}  For essentially all quadratic $f \in \Z[x]$, either $f^n(\gamma) > 0$ for all $n$ large enough or $f^n(\gamma) < 0$ for all $n$ large enough.  In this case, the proof of Theorem \ref{mainone}, part 1, shows that for $f$ satisfying the hypotheses of the theorem, $G_n(f)$ has finite index in the full automorphism group of the complete infinite binary rooted tree.

\smallskip

{\bf \noindent Remark:}  With only minor modifications, the proof of Theorem \ref{mainone} applies to the case of $f \in \calO_K[x]$, where $\calO_K$ is the ring of integers in a number field $K$.  In this case, there can be many more $f$ that satisfy the hypotheses of the theorem.  For instance, if we let $\alpha$ be a root of
$x^3 + 3x^2 + 2x + 1$ and $K = \Q(\alpha)$, then not only do the families in Theorem \ref{maintwo} have $\{f^n(0)\}$ finite and not containing $0$ (taking $k \in \calO_K$ rather than $k \in \Z$), but so does $f(x) = x^2 + x + \alpha$.  Indeed, we have $\{f^n(0)\} = \{\alpha, \alpha^2 + 2\alpha, -1\}$.  \label{remkmainone}

\section{Irreducibility and Stability of Polynomial Iterates} \label{stable}

In order to deduce Theorem \ref{maintwo} from Theorem \ref{mainone}, we must show that for $f$ belonging to the four families in question, $f$ is critically infinite and $g \circ f^n$ is irreducible for all $n$ and some appropriate $g$.  This section examines the irreducibility question, which has been studied in its own right, e.g. in \cite{fein}, \cite{ayad}.

\begin{definition} Let $K$ be a field and $f$ and $g$ be polynomials in $K[x]$.  We say $g$ is {\em $f$-stable} if $g \circ f^n$ is irreducible over $K$ for $n = 0,1,2, \ldots$ (note that in particular $g$ is irreducible over $K$).  We say $f$ is {\em stable} if $f$ is $f$-stable, i.e., all iterates of $f$ are irreducible.  We say $f$ is {\em eventually stable} if some iterate of $f$ is a product of $f$-stable polynomials.
\end{definition}

\begin{lemma}[Capelli's Lemma] Let $K$ be a field, $f(x), g(x) \in K[x]$, and let $\beta \in \overline{K}$ be any root of $g(x)$.  Then $g(f(x))$ is irreducible over $K$ if and only if both $g$ is irreducible over $K$ and $f(x) - \beta$ is irreducible over $K(\beta)$.
\end{lemma}

For a proof, see \cite[Lemma 0.1]{fein}.  The following proposition is a remarkably useful tool for determining $f$-stability when $f$ is quadratic.

\begin{proposition} \label{stabcrit}
Let $f(x) = ax^2 + bx + c \in K[x]$ and let $\gamma = -b/2a$ be the critical point of $f$.  Suppose that $g \in K[x]$ is such that $g \circ f^{n-1}$ is irreducible over $K$ for some $n \geq 2$.
Then $g \circ f^n$ is irreducible over $K$ if $g(f^n(\gamma))$ is not a square in $K$.
\end{proposition}

\Pf (cf \cite[Lemma~4.13]{galmart})  Denote by $d$ the degree of $g \circ f^{n-1}$.  By Capelli's Lemma and the irreducibility of $g \circ f^{n-1}$, we have $g \circ f^n$ irreducible if for any root $\beta$ of $g \circ f^{n-1}$, $\Disc f(x) - \beta = b^2 - 4ac + 4a\beta$ is not a square in $K(\beta)$.  This must hold if $N_{K(\beta)/K} (b^2 - 4ac + 4a\beta)$ is not a square in $K$.  But
\begin{eqnarray*}
N_{K(\alpha)/K} (b^2 - 4ac + 4a\beta) & = & (-4a)^{d} \prod_{\text{$\beta$ root of $g \circ f^{n-1}$}}
\left(- \frac{b^2}{4a} + c \right) - \beta \\
& = & (-4a)^{d} \cdot g(f^{n-1}(-b^2/4a + c)) = (-4a)^{d} \cdot g(f^{n-1}(f(\gamma))).
\end{eqnarray*}
Since $n \geq 2$, we have $d$ even, which proves the proposition.  \done

From Proposition \ref{stabcrit}, it follows that stability is a generic property for quadratic polynomials with integral coefficients.  As an example of the kind of result that Proposition \ref{stabcrit} allows, we give the following lemma and theorem, which show that for quadratic polynomials with $f^2$ irreducible and given critical point $\gamma$, all but finitely many are stable.

\begin{lemma} \label{long}
Let $f \in \Z[x]$ be monic, quadratic, and nonsquare, and write $f = (x - \gamma)^2 + \gamma + m$.  Suppose that
$\abs{m} > 6 + 3\sqrt{\abs{\gamma} + 1}$.  Then for every $n \geq 3$, $f^n(\gamma)$ is not a square in $\Q$.  Furthermore, if $\gamma \in \Z$ then the same conclusion holds if $\abs{m} > 1 + \sqrt{\abs{\gamma}+1}$.
\end{lemma}

\Pf  Note that either 1) $\gamma \in \Z$ or 2) $\gamma \not\in \Z$ but $2\gamma \in \Z$.  In either case we have by hypothesis $\abs{m} > 2$, and since the denominator of $m$ cannot exceed $4$ we in fact have $\abs{m} \geq 9/4$.    For $r/s \in \Q$ we have
$$f \left(\frac{r}{s}\right) = \frac{1}{s^2} \left[ (r-s\gamma)^2 + s^2(\gamma + m) \right],$$
and if we assume that $s \gamma$ and $s^2(\gamma + m)$ are integers then
$f(r/s)$ is not a square if
\begin{equation} \label{monk2}
\abs{s^2(\gamma + m)} < 2 \abs{r - s\gamma} - 1.
\end{equation}
This follows because $\gamma + m \neq 0$ (since $f$ nonsquare) and
$$\min{\left\{\abs{x^2 - (x-1)^2}, \abs{x^2 -(x + 1)^2}\right\}} \geq 2\abs{x} - 1$$
for $x \in \Z$.

Put $\frac{r}{s} = f^{n-1}(\gamma)$, and note that
$s = j^{2^{n-1}}$, where $j$ is 1 if $\gamma \in \Z$ and $2$ otherwise.  Dividing \eqref{monk2} by $s^2$ then yields that $f^n(\gamma)$ is not a square if
\begin{equation} \label{prince}
\abs{\gamma + m} < j^{-2^{n-1}} \left( 2 \abs{f^{n-1}(\gamma) - \gamma} - 1 \right).
\end{equation}
A quick induction shows that $f^{n-1}(\gamma) - \gamma = f_m^{n-1}(0)$, where
$f_m(x) = x^2 + m$.

We now show that \eqref{prince} holds for $n \geq 3$.  First we show that it holds for $n = 3$; in this case
\eqref{prince} becomes
\begin{equation} \label{prince2}
j^{4}\abs{\gamma + m} <  2 \abs{m^2 + m} - 1.
\end{equation}
Suppose $j = 2$, or equivalently that $\gamma \not\in \Z$.  We have $\abs{m} - 6 >  3\sqrt{\abs{\gamma}+1}$, whence by squaring both sides and subtracting we have $\abs{m}^2 - \abs{m} > 11\abs{m} + 9 \abs{\gamma} + 9 - 36$.  Since $\gamma \not\in \Z$, we have $\gamma \geq 1/2$.  Thus $\abs{m} \geq 9 + 3/4$ (since the denominator of $m$ is $4$), and we have
$11\abs{m} + 9 \abs{\gamma} -27 > 8 \abs{m} + 8\abs {\gamma} + 1$.
Now
$$2 \abs{m^2 + m} - 1 \geq 2(\abs{m}^2 - \abs{m}) - 1 > 16\abs{m} + 16\abs{\gamma} + 1 >
16\abs{\gamma + m},$$
and thus \eqref{prince2} is satisfied.
In the case $j=1$, we have $(\abs{m}-1)^2 > \abs{\gamma} + 1$, and
thus $\abs{m}^2 - \abs{m} > \abs{\gamma} + \abs{m}.$
Therefore
$$2\abs{m^2+m} - 1 \geq 2(\abs{m}^2 - \abs{m}) -1 > 2(\abs{\gamma} + \abs{m}) -1 > \abs{\gamma} + \abs{m},$$
where the last inequality follows since $\abs{m} \geq 9/4$.   Hence \eqref{prince2} is satisfied.

To show that \eqref{prince} holds for $n > 3$, we must show the function
$j^{-2^n} \left( 2 \abs{f_m^n(0)} - 1 \right)$ increases with $n$, i.e.
\begin{equation}
j^{2^{n-1}} < \frac{2 \abs{f_m^n(0)} - 1}{2 \abs{f_m^{n-1}(0)} - 1}
\end{equation}
Noting that $f^n_m(0) = (f_m^{n-1}(0))^2 + m$, it is enough to prove
$$j^{2^{n-1}} <  \abs{f_m^{n-1}(0)} - \frac{2\abs{m} + 1}{2 \abs{f_m^{n-1}(0)}}.$$
Since $\abs{m} \geq 9/4$, it follows that $2\abs{m} + 1 < 2 \abs{f_m^{n-1}(0)}$ for all $n \geq 3$.
Thus it suffices to show
\begin{equation} \label{thingtoshow}
\abs{f_m^{n-1}(0)}  > 1 + j^{2^{n-1}}
\end{equation}
for all $n \geq 4$.  One can easily use $\abs{m} > 9$ to show that this holds for $n = 4$ (in the case $j = 1$, $\abs{m} > 2$ suffices).   Now assume that \eqref{thingtoshow} holds for $n-1$, with $n \geq 5$.  We show that it holds for $n$.

If $j=1$, then $\abs{m} > 2$ guarantees that $\abs{f_m^n(0)} > \abs{f_m^{n-1}(0)}$, and we are done.
If $j=2$, suppose first that $\abs{m} > 2j^{2^{n-1}}$.   Since $\abs{m} > 3$ (indeed we know $\abs{m} > 9$), we have $\abs{m}^2 - 3\abs{m} > 0$, and this implies that $\abs{m^2} - \abs{m} > 2\abs{m} > \abs{m} + 1 + j^{2^{n-1}}$.
Therefore $\abs{f_m^2(0)} > \abs{m} + 1 + j^{2^{n-1}}$.  Squaring both sides leads to
$\abs{f_m^3(0)} > 1 + j^{2^{n}}$.  Since $\abs{f_m^n(0)} > \abs{f_m^{n-1}(0)}$, we have
$\abs{f_m^n(0)} > 1 + j^{2^{n}}$.

Now suppose that $\abs{m} < 2j^{2^{n-1}}$ (equality is impossible since $j = 2$ implies $\gamma \not\in \Z$).  Given $\abs{f_m^{n-1}(0)}  > 1 + j^{2^{n-1}}$, squaring both sides gives
$$\abs{f_m^n(0) - m} > 1 + 2j^{2^{n-1}} + j^{2^{n}}  >  1 + j^{2^{n}} + \abs{m}.$$
Therefore $\abs{f_m^n(0)} + \abs{m} >  1 + j^{2^{n}} + \abs{m},$
and we are done.  \done

\begin{theorem} \label{stabthm}
Let $f \in \Z[x]$ be monic, quadratic, and irreducible, and write
$f = (x-\gamma)^2 + \gamma + m.$
Suppose that $\abs{m} > 6 + 3\sqrt{\abs{\gamma} + 1}$ (if $\gamma \in \Z$ then $\abs{m} > 1 + \sqrt{\abs{\gamma} + 1}$ suffices), and that
\begin{equation} \label{f2irr}
\frac{-m \pm \sqrt{f^2(\gamma)}}{2} \not\in \Q^{*2}.
\end{equation}
Then $f$ is stable.
\end{theorem}

{\bf \noindent Remark}: Condition \eqref{f2irr} is always true if $\gamma \not\in \Z$, as $f^2(\gamma)$ is an odd number divided by 16.

\Pf Given that $f$ is irreducible, condition \eqref{f2irr} is equivalent to $f^2$ being irreducible, as a quick calculation and application of Proposition \ref{stabcrit} shows.  Since $f$ is irreducible, it is not a square, whence Lemma \ref{long} shows that $f^n(\gamma)$ is not a square for all $n \geq 3$.  By induction and Proposition \ref{stabcrit} this proves that $f$ is stable.   \done

\begin{proposition} \label{famonestab}
Let $k \in \Z \backslash \{0,-1\}$, and put $f = x^2 + k$.  Then $f$ is stable if $-k$ is not a square in $\Z$ and $f$ factors as the product of two linear $f$-stable polynomials if $-k$ is a square in $\Z$.
\end{proposition}

\Pf Suppose first that $-k$ is not a square in $\Z$, whence $f$ is irreducible.  The critical point of $f$ is $\gamma = 0$, and $f^2(0) = k^2+k$, which is not a square since $k \not\in \{0,-1\}$.  Thus the condition \eqref{f2irr} is satisfied, as indeed are all the hypotheses of Theorem \ref{stabthm} as long as
$k \neq 1$, $k \neq 2$, and $k \neq -2$.  For these cases, we verify \eqref{prince} directly for all $n \geq 2$.  If $k=1$, then \eqref{prince} becomes $1 < 2\abs{f^n(0)} - 1$, which clearly holds for $n \geq 2$.  A similar argument applies to $k = 2$.  If $k = -2$ the forward orbit of $0$ lies in $\{\pm2\}$, whence $f^n(0)$ is not a square for $n \geq 2$.  In all three cases induction and Proposition \ref{stabcrit} show $f$ is  stable.

If $-k = a^2$ for some $a \in \Z$ then we have $f = g_1g_2$, where $g_1 = x+a$ and $g_2 = x-a$.
Note that since $k \not\in \{0,-1\}$, we have $\abs{a} \geq 2$.  To show that
$g_i(f^{n}(0))$ is not a square for all $n \geq 2$, we note $g_i(f^{n}(0)) = f^{n-1}(0)^2 - a^2 \pm a$, so by the same reasoning used to obtain \eqref{prince} we have
$g_i(f^{n}(0))$ not a square if $a^2 \pm a < 2 \abs{f^{n-1}(0)} - 1$.  Since $\abs{a} \geq 2$ this holds for $n = 2$.  Moreover, again because $\abs{a} \geq 2$, it follows that
$\abs{f^{n+1}(0)} > \abs{f^{n}(0)}$ for all $n \geq 2$.  Induction and Proposition \ref{stabcrit} now show $g_1$ and $g_2$ are $f$-stable.  \done

\begin{proposition} \label{famtwostab}
Let $k \in \Z \backslash \{0\}$, and put $f = x^2 + kx - k$.  Then $f$ is stable if $k \neq 4$ and $f$ is the square of an $f$-stable polynomial if $k=4$.
\end{proposition}

\Pf  In the notation of Theorem \ref{stabthm}, we have $\gamma = -k/2$ and $m = k(2-k)/4 = \gamma - \gamma^2$.  Thus $f$ is irreducible as long as $2\gamma - \gamma^2 \neq 0$, which is the case when $k \neq 0, 4$.  To show that $f^2$ is irreducible, we note that $f^2(\gamma) =
(\gamma - \gamma^2)^2 + 2\gamma - \gamma^2$, so to ensure $f^2(\gamma)$ is not a square, it is enough to show $\abs{2\gamma - \gamma^2} < 2\abs{\gamma - \gamma^2} - 1$.  This holds for
$\abs{\gamma} > 2$.  Now to use Theorem \ref{stabthm} we must find the $\gamma$ with
\begin{equation} \label{thingy}
\abs{\gamma - \gamma^2} > 6 + 3 \sqrt{\abs{\gamma}+1}.
\end{equation}
A derivative shows that $\abs{\gamma - \gamma^2} - 6  - 3\sqrt{\abs{\gamma}+1}$ increases with $\abs{\gamma}$ provided $\abs{\gamma} > 2$.  One easily sees that \eqref{thingy} holds for all $\gamma \not\in (-7/2, 9/2)$.  Similarly, a derivative shows that
$\abs{\gamma - \gamma^2} - 1  - \sqrt{\abs{\gamma}+1}$ increases with
$\abs{\gamma}$ provided $\abs{\gamma} > 1$.  A quick calculation now shows
$\abs{\gamma - \gamma^2}  > 1 + \sqrt{\abs{\gamma}+1}$ for all $\gamma \not\in (-2, 3)$.
Therefore by Theorem \ref{stabthm} we have that $f$ is stable unless $k = -5, -3, -2, -1, 1, 2, 3, 5, 7$.
For all of these save $k = \pm 1$ one can easily check that all iterates of $f$ are Eisenstein and thus irreducible.  The case $k = -1$ is proven in \cite[Proposition 1.1]{odoniwn}.  The case $k = 1$ is handled in Proposition \ref{famthreestab}.

For $k = 4$, we clearly have $f = g^2$, with $g = (x-2)$.  Moreover, $f^n(\gamma) = 0$ for $n = 1$ and
$f^n(\gamma) = 4$ for $n \geq 2$.  Thus $g \circ f^n(\gamma)$ is not a square for all $n \geq 1$, and it follows that $g \circ f^n(\gamma)$ is irreducible for all $n \geq 1$.  \done

\begin{proposition} \label{famthreestab}
Let $k \in \Z \backslash \{0\}$, and put $f = x^2 + kx - 1$.  Then $f$ is stable if $k \neq -1$ and
$f^3$ is the product of two $f$-stable polynomials if $k = -1$.
\end{proposition}

\Pf In the notation of Theorem \ref{stabthm}, we have $\gamma = -k/2$ and $m = (-k^2 + 2k - 4)/4$.
We have that $f$ is irreducible as long as $k^2 + 4$ is not a square, which holds for all $k \neq 0$.
To show $f^2(x)$ is irreducible, note that
$$f^2(\gamma) = \frac{1}{16}(k^4 - 4k^3 + 8k^2 + 16k) =
\left( \frac{1}{4}(k-1)^2 \right)^2 + \frac{1}{16}(2k^2 - 12k -1).$$
Thus $f^2(\gamma)$ is not a square as long as $\abs{2k^2 - 12k -1} < 8 \abs{k-1}^2 - 1$, which holds for
$k \not\in \{0, 1, 2, 3\}$.

One checks that $\abs{(-k^2 + 2k - 4)/4} > 6 + 3\sqrt{\abs{k/2} + 1}$ as long as $k \not\in (-7, 9)$ and also
$\abs{(-k^2 + 2k - 4)/4} > 1 + \sqrt{\abs{k/2} + 1}$ as long as $k \not\in (-2,4)$.
Using Theorem \ref{stabthm}, we now have that $f$ is stable as long as $k \not\in \{-5,-3,-1,0,1,2,3,5,7\}$.
Consider $k = 7$.  We have $\gamma = -7/2$, $m = -39/4$.  We note that in the proof of Lemma \ref{long} we need only verify \eqref{prince} for some $n$ and \eqref{thingtoshow} for all $i \geq n$.  When $n = 4$ we have $f^{n-1}(\gamma) - \gamma > 7268$ and $2^{2^{n-1}} \cdot \abs{\gamma + m} = 3392$, whence \eqref{prince} holds.  To verify \eqref{thingtoshow} for $n \geq 4$, note that it holds for $n=4$ since $f_m^3(0) > 7268 > 1 + 2^8$.  Moreover, given that \eqref{thingtoshow} holds for $n$, to show it holds for $n+1$ it is enough to have $\abs{m} > 3$ (see proof of Lemma \ref{long}).  Hence we have shown that $f^n(\gamma)$ not a square for $n \geq 4$.  One checks the remaining cases by hand, proving that $f$ is stable.  In a similar manner, one verifies that $f$ is stable for $k \in \{-5, -3, 5, 7\}$, although for $k = -3$ one must take $n = 6$ and for $k = 5$ one must take $n = 7$.

For $k \in \{1, 2, 3\}$, let $\alpha_1 < \alpha_2$ be the two (real) roots of $f$, and put
$I = (\alpha_1, 0)$.  Using the fact that $f(\gamma) \in I$, it follows that $f(I) \subseteq I$, whence
$f^n(\gamma) < 0$ for all $n$.  Thus $f^n(\gamma)$ cannot be a square, so $f$ is stable.  This leaves the case $k = -1$, where we see that although $f$ and $f^2$ are irreducible,
$$f^3(x) = (x^4 - 3x^3 + 4x - 1)(x^4 - x^3 - 3x^2 + x + 1) := g_1(x)g_2(x),$$
where $g_1$ and $g_2$ are irreducible.  To prove that $g_1$ is $f$-stable, it is enough to show that
$g_1(f^n(\gamma))$ is not a square for all $n \geq 1$.  Note that $g_1(-1) = -1$ and also $f^2(-1) = -1$.  Working modulo 3, we have $\gamma \equiv -1$, whence $g_1(f^n(\gamma)) \equiv -1 \bmod{3}$ for $n$ even.  Working modulo 19, we have $f^5(\gamma) \equiv -1 \bmod{19}$, whence
$g_1(f^n(\gamma)) \equiv -1 \bmod{19}$ for $n$ odd, $n \geq 5$.  One now checks that $g_1(f^n(\gamma))$ is not a square in $\Q$ for $n = 1, 3$, showing that $g_1$ is $f$-stable.  For $g_2$ we have
$g_2(1) = g_2(-1) = -1$, and since $f^n(\gamma) \equiv \pm 1 \bmod{3}$ for all $n$, we have
$g_2(f^n(\gamma)) \equiv -1 \bmod{3}$ for all $n$, proving that $g_2$ is $f$-stable.  \done

We note that it is possible for a quadratic $f$ to be irreducible over $\Q$ yet for $f^2$ to have a nontrivial factorization.  One example is $x^2 + 10x + 17$, which is irreducible yet clearly fails to satisfy condition $\eqref{f2irr}$.  Moreover, as in the case of $f(x) = x^2 - x - 1$, it is also possible for $f$ and $f^2$ to be irreducible yet for $f^3$ to have a nontrivial factorization.  Computer searches have revealed no monic, quadratic $f \in \Z[x]$ whose first three iterates are irreducible yet whose fourth is not.  Moreover, in all cases checked, the reducible iterate has factored into two $f$-stable polynomials.  Generalizing this observation, we make the following conjecture:

\begin{conjecture} \label{evstabconj}
Let $f \in \Z[x]$ be monic and quadratic, and suppose that $0$ is not periodic under $f$.  Then $f$ is eventually stable.
\end{conjecture}

We remark that if $0$ is periodic under $f \in \Z[x]$, then $f$ cannot be eventually stable.  This is because the factors of $f^n$ are the minimal polynomials over $\Q$ of the preimages of $0$
(in $\overline{\Q}$) under $f$.  Thus if $0$ is periodic under $f$, $x$ occurs as a factor of some $f^m$, and it follows that $f^n$ has a linear factor for each $n \geq 1$.  Thus for instance no member of the family $x^2 + kx - (k+1), k \in \Z$, is eventually stable.

\section{results for specific families} \label{families}
In this section we draw together the results from all the previous sections to prove Theorem \ref{maintwo}, treating one family at a time.  We also give Conjecture \ref{bigconj}, which treats more general $f$ than those meeting the hypotheses of Theorem \ref{mainone}.

\begin{theorem} \label{main1}
Let $f = x^2 -kx + k$ for $k \in \Z$.  Then $D(P(f,a_0)) = 0$ for all $a_0$.
\end{theorem}

\Pf Clearly $P(f,a_0)$ is finite for $k=0$, so this case is covered.  To apply Theorem \ref{mainone} we must show that $f$ is stable and critically infinite.  If $k \neq 4$, then by Proposition \ref{famtwostab} we have that $f$ is stable.  We now show that $\{f^n(\gamma) : n \geq 1\}$ is infinite for $k \not\in \{-2, 2, 4\}$. Recall that we denote the critical point of $f$ by $\gamma$, which is $k/2$ for the family under consideration.  Suppose first that $k$ is odd.  Letting $a_{n-1}$ and $a_n$ denote the numerators of $f^{n-1}(\gamma)$ and $f^n(\gamma)$ respectively, we have
$a_n = a_{n-1}^2 - k2^{2^{n-1}}a_{n-1} + k2^{2^n}$.  Moreover, $f^n(\gamma) = a_n/2^{2^n}$.  Thus by induction $a_n$ must be odd, and therefore $f$ is critically infinite since the denominator of $f^n(\gamma)$ in lowest terms is increasing.

Now suppose that $k$ is even. We note that $f$ is increasing on $(\gamma, \infty)$, and that the largest fixed point of $f$ is  $\rho = \max\{1, k\}$.  If $f^n(\gamma) > \rho$ for some $n$, one easily sees that
$f^{m+1}(\gamma) > f^m(\gamma)$ for all $m \geq n$, showing that $f$ is critically infinite.
We show $f^n(\gamma) > \rho$ for $n=2$.  Note that $f^2(\gamma) = \frac{1}{16}(k^4-4k^3) + k$.  If $k > 0$ then $\rho = k$, and we have $f^2(\gamma) > \rho$ as long as $k > 4$.  If $k < 0$ then
$\rho = 1$, and we have $f^2(\gamma) > \rho$ as long as $k^4 - 4k^3 + 16k - 16 > 0$.  This is equivalent to $(k-2)^3(k+2) > 0$, which holds when $k < -2$.  Thus $f$ is critically infinite as long as
$k \notin \{-2, 2, 4\}$, and the theorem is proved in all but these cases.

For $k \in \{-2, 2, 4\}$, we first show that $\#H_n(f,g) = 2$ for all $n$.
If $f = x^2-4x+4$ and $g = x-2$ then roots of $g \circ f^n$ are the same as roots of iterates of $x^2 - 2$,
that is, $\zeta_{2^{n+2}} + \zeta_{2^{n+2}}^{-1}$, and the extensions thus generated are $\Q(\zeta_{2^{n+2}}) \cap \mathbb{R}$.  In the remaining two cases we have $g = x$. If $f=x^2 + 2x - 2$, iterates of $f$ form extensions that are $\Q(\zeta_{3 \cdot 2^{n+1}}) \cap \mathbb{R}$. Finally, roots of iterates of $x^2 -2x + 2$ are roots of $x^{2^n} + 1$ for all $n \geq 1$.  In all three cases one easily shows that $\#H_n(f,g) = 2$ for all $n$.
We now show by induction that the identity is the only element of $G_n(f,g)$ that fixes a root of $g \circ f^n$, whence the proportion of such elements is $2^{-n}$ and from Theorem \ref{upbound} we have $D(P(f, a_0)) = 0$.  The case $n=1$ is clear.  Assuming that only the identity in $G_{n-1}(f,g)$ fixes any roots of $g \circ f^{n-1}$, we have that only elements of $H_n(f,g)$ can fix any roots of $g \circ f^n$.  Since $H_n(f,g)$ is not trivial, $g \circ f^n$ factors over
$\Q(g \circ f^{n-1})$ as a product of quadratics, and $H_n(f,g)$ must act transitively on the roots of each of these factors.  Thus $H_n(f,g)$ consists of the identity and an element that exchanges the roots of each quadratic factor of $g \circ f^n$ over $\Q(g \circ f^{n-1})$, and thus fixes no roots of $g \circ f^n$.
\done

\begin{theorem} \label{main2}
Let $f = x^2 +kx - 1$ for $k \in \Z \backslash \{0,2\}$.  Then $D(P(f,a_0)) = 0$ for all $a_0$.
\end{theorem}

\Pf Part (1) of Theorem \ref{mainone} applies immediately (with $g(x) = x$) in the case that $f$ is stable and critically infinite.  If $k \neq -1$, then by Proposition \ref{famthreestab} we have that $f$ is stable.  We now show that $\{f^n(\gamma) : n \geq 1\}$ is infinite for $k \not\in \{0, 2\}$.  As in the proof of Theorem \ref{main1}, the case $k$ odd is covered.

Now suppose that $k$ is even, and note that $\gamma = -k/2$, $m = (-k^2 + 2k - 4)/4$.  Following the proof of Theorem \ref{main1}, we have here that the largest fixed point $\rho$ of $f$ is at most
$\max\{1, -k+2\}$.  Again it is enough to show $f^n(\gamma) > \rho$ for some $n$, and we take $n = 2$.
We have $f^2(\gamma) = \frac{1}{16}(k^4-4k^3 + 8k^2 - 16k)$, which is greater than 1 for $k \geq 4$ and greater than $-k + 2$ for $k \leq -2$.

This leaves the case $k = -1$.  Since from Proposition \ref{famthreestab} we have $f^3 = g_1g_2$ with $g_1$ and $g_2$ $f$-stable, Part (1) of Theorem \ref{mainone} applies to show that both $P(g_1,f, a_0)$ and $P(g_2,f, a_0)$ have density $0$ for every $a_0 \in \Z$.  Since $P(f, a_0)$ is the union of these two sets, it too has density $0$.  The theorem is proved.  \done

The cases of $x^2 - 1$ and $x^2 + 2x - 1 = (x+1)^2 - 2$ present peculiar difficulties, as both are critically finite, but the relevant Galois groups do not lend themselves to easy analysis.  We remark that to show $D(P(x^2-1, a_0)) = 0$, put $f = x^2 - 1$ and note that $f^n(a_0) \equiv 0 \bmod{p}$ for some $n$ implies either 1) $a_0 \equiv 0$ or $ - 1 \bmod{p}$ or 2) $f^m(a_0) \equiv 1 \bmod{p}$ for some $m$.  Thus it suffices to show that the density of primes $p$ with $f^{n}(a_0) \equiv 1 \bmod{p}$ is zero, and this may be accomplished by considering the Galois groups of $f^n(x) - 1$.  These are the same as those generated by iterates of $h(x) = (x+1)^2 - 2$.  Iterates of $h$ generate extensions that are finitely ramified (indeed, unramified outside $\{2\}$), yet are much larger than the extensions generated in the exceptional cases handled in Theorem \ref{main1}.  There appears to be no easy way to show that $\inflim{n} \P(X_n) = 0$ in this case.

We now turn to the case where $f$ is monic and quadratic and all primes dividing an element of the critical orbit are recurrent.  This is the case, for instance, whenever $0$ is the critical point, as the following lemma shows.  See \cite[Section 3]{brian} for a nice discussion of matters similar to the following lemma and proposition.

\begin{lemma} \label{rigidlem} Let $f = x^2 + k$ for $k \in \Z \backslash \{0,-1\}$, and put $f^n(0) = t_n$ for $n \geq 1$.  Then $t_n$ is a rigid divisibility sequence.
\end{lemma}

\Pf (cf \cite[Lemma 6.1]{galmart})
Since we have excluded $k \in \{0,-1\}$, we have $t_n \neq 0$ for all $n$.
We use induction on $j$ to show $v_p(t_n) = e > 0$ implies
$v_p(t_{nj}) = e$ for all $j \geq 1$.  The case $j=1$ is trivial.  Suppose
inductively that ${v}_{p}(t_{n(j-1)}) = e$.  Note that
$t_{nj} = f^{n(j-1)}(t_n)$, and also $f^{n(j-1)}$ is a polynomial in $x^2$.
Thus we can write
$$f^{n(j-1)}(x) = x^2g(x) + f^{n(j-1)}(0) = x^2g(x) + t_{n(j-1)},$$
for some $g \in \Z[x]$.  Hence putting $x = t_n$ we have
$t_{nj} = t_n^2g(t_n) + t_{n(j-1)}.$
Now ${v}_{p}\left[(t_n)^2(g(t_n))\right] \geq 2e$, and by our inductive hypothesis
${v}_{p}(t_{n(j-1)}) = e$.  Since $e \geq 1$, the first summand vanishes to
higher order at $p$ than the second, so we conclude ${v}_{p}(p_{nj}) = e$.  \done

If $f^n(0)$ is a rigid divisibility sequence and $f$ is reducible, its factors give rise to other rigid divisibility sequences, as the following Proposition shows.

\begin{proposition} \label{rigid}
Suppose that $f \in \Z[x]$ has all iterates separable, and that $0$ is a critical point of $f$ and $(f^n(0))$ is a rigid divisibility sequence with infinitely many values.  If $g$ divides $f$
then the sequence $g(0), g(f(0)), g(f^2(0)), \ldots$ is a rigid divisibility sequence.
\end{proposition}

\Pf Suppose that $gh = f$, and for $n \geq 1$ put $a_n = f^n(0)$ and $b_n, c_n$ respectively for
$g(f^{n-1}(0)), h(f^{n-1}(0))$.  Recall that by convention $f^0(x) = x$, and note that we have
$b_nc_n = a_n$ for all $n \geq 1$.  If $v_p(b_i) = e > 0$, then we have $a_i = f^i(0) \equiv 0 \bmod{p^e}$, whence $0$ is periodic $\bmod{p^e}$ under $f$, with period dividing $i$.  Hence in particular
$f^{mi-1}(0) \equiv f^{i-1}(0) \bmod{p^e}$ for all $m \geq 1$.  Now $v_p(b_i) = e > 0$ implies that $g(f^{i-1}(0)) \equiv 0 \bmod{p^e}$, and thus $f^{i-1}(0)$ is a root of $g \bmod{p^e}$.  Therefore
$f^{mi-1}(0)$ is also a root of $g \bmod{p^e}$, implying that $v_p(b_{mi}) \geq e$.  Similarly, if
$v_p(c_i) = e > 0$, then $v_p(c_{mi}) \geq e$ for all $m \geq 1$.

Now suppose that $v_p(b_i) = e_1 > 0$.  Then $v_p(c_i) = e_2 \geq 0$, and $v_p(a_i) = e_1 + e_2$.
Since $(a_n)$ is a rigid divisibility sequence, we have $v_p(a_{mi}) = e_1 + e_2$ for all $m \geq 1$.
But by the previous paragraph $v_p(b_{mi}) \geq e_1$ and $v_p(c_{mi}) \geq e_2$.  Thus
$a_{mi} = b_{mi}c_{mi}$ implies $v_p(b_{mi}) = e_1$ and $v_p(c_{mi}) = e_2$.  Hence $(b_n)$ is a rigid divisibility sequence, as desired. \done

\begin{theorem} \label{main3}
Let $f = x^2 + k \in \Z[x]$ with $k \neq -1$.  Then $D(P(f,a_0)) = 0$ for all $a_0 \in \Z$.
\end{theorem}

\Pf First suppose that $-k$ is not a square in $\Z$.
By Proposition \ref{famonestab} $f$ is stable, and provided that $k \neq -2$ we have $f$ critically infinite, whence the theorem follows from Theorem \ref{mainone}, part 2, taking $g(x) = x$.
The case $k = -2$ may be resolved directly, as has already been done in the proof of Theorem \ref{main1}.

Suppose now that $-k = c^2$ for some $c \in \Z$ with $\abs{c} \geq 2$, thereby excluding $k \in \{0, -1\}$. Then $\{f^n(0)\}$ is infinite.  Letting $g(x) = x - c$, we have that $g \circ f^n$ is irreducible for all $n \geq 1$ by Proposition \ref{famonestab}.  Moreover, the critical point of $f$ is $\gamma = 0$ and Lemma \ref{rigidlem} and Proposition \ref{rigid} imply that $g(f^n(0))$ is a rigid divisibility sequence.   Thus part 2 of Theorem \ref{main2} shows that $D(P(g,f, a_0)) = 0$ for all $a_0 \in \Z$.  Similarly, $D(P(h,f, a_0)) = 0$, where $h(x) = x + c$.  Since $g(f^n(a_0)) \cdot h(f^n(a_0)) = f^{n+1}(a_0)$, we have that $P(f, a_0)$ is contained in the union of two density zero sets, and thus has density zero.

Finally, if $k = 0$ it is clear that the set of divisors of $a_n$ is finite.  \done

For our final family, we ride the coattails of Theorem \ref{main3}.  This is very much along the lines of the sequence factorization described in \cite[Section 3]{brian}.

\begin{theorem} \label{main4}
Let $f = x^2 - 2kx + k \in \Z[x]$ with $k \neq \pm 1$.  Then $D(P(f,a_0)) = 0$ for all $a_0 \in \Z$.
\end{theorem}

\Pf Note first that $f(x) = (x - k)^2 - k^2 + k$, so that $f = \phi \circ h \circ \phi^{-1}$, where
$\phi(x) = x + k$ and $h(x) = x^2 - k^2$.  Hence $f^n = \phi \circ h^n \circ \phi^{-1}$, implying that
$f^n(a_0) = h^n(a_0 - k) + k$.  But $(h^n(a_0 - k) + k)(h^n(a_0 - k) - k) = h^{n+1}(a_0 - k)$, and
$D(P(h, a_0 - k)) = 0$ by Theorem \ref{main3} provided $k^2 \neq 1$.  \done

We have now completed the proof of Theorem \ref{maintwo}, in which Theorem \ref{mainone} has been the driving force.  Unfortunately, quadratic $f \in \Z[x]$ meeting the hypotheses of Theorem \ref{mainone} are rare (and indeed remain rare despite growing more numerous if $\Z$ is replaced by the ring of integers of a number field).  However, we conjecture that the conclusion of Theorem \ref{mainone} holds for generic quadratic $f \in \Z[x]$:

\begin{conjecture} \label{bigconj}
Let $f \in \Z[x]$ be quadratic, and suppose that $f$ is stable and critically infinite.  Then
$D(P(f, a_0)) = 0$ for all $a_0 \in \Z$.
\end{conjecture}

The obstacle impeding a proof of Conjecture \ref{bigconj} is a knowledge of the arithmetic of the critical orbit of $f$.  It is enough to show that given a stable, quadratic $f \in \Z[x]$ and $b \in \frac{1}{2}\Z$ with $\{f^n(b) : n \geq 1\}$ infinite, then there are infinitely many $n$ such that there exists a prime $p$ with $v_p(f^n(b))$ odd and $v_p(f^m(b)) = 0$ for all $m < n$.  This assertion is plausible in light of work of \cite{everest} showing that for any nonzero $c \in \Z$, ``most'' $n$ have the property that there exists $p$ with $v_p(n^2 + c) = 1$ but $v_p(m^2 + c) = 0$ for all $m < n$.  It seems plausible as well that a similar assertion is true for translated sequences $\{g(f^n(b)) : n \geq 1\}$.  This would prove that Conjecture \ref{bigconj} holds for all quadratic $f$ that are eventually stable and critically infinite.  If Conjecture \ref{evstabconj} is true as well, this would show that $D(P(f, a_0)) = 0$ for virtually all quadratic $f$.

\section{The infinitude of prime divisors of first-order quadratic recurrences} \label{inf}

\begin{theorem}
Suppose that $f \in \Z[x]$ has degree 2, $f \neq x^2$, and $a_0 \in \Z$ is such that $\{f^n(a_0) : n \in \N\}$ is infinite.  Then $P(f, a_0)$ is infinite.
\end{theorem}

\Pf First suppose that $\Disc f \neq 0$.  Let $p_1, \ldots, p_r$ be primes, and recall $a_n = f^n(a_0)$.  We show that there can be only finitely many
$n$ such that $a_n$ is a product of powers of the $p_i$.  Write $f = b_2x^2 + b_1x + b_0$, and suppose
$a_n = p_1^{m_1} \cdots p_r^{m_r}$ for some $m_1, \ldots m_r$.  Then there is a $c_n$ such that $a_n = c_n^3 p_1^{j_1} \cdots p_r^{j_r}$, with $0 \leq j_i \leq 2$ for each $i$.  Put $d(a_n) =  p_1^{j_1} \cdots p_r^{j_r}$.  We thus have that $(a_{n-1}, c_n)$ gives a solution to
$$C : b_2x^2 + b_1x + b_0 = d(a_n)y^3.$$
By completing the square on the left-hand side, and using that $f$ is not already a square since it has nonzero discriminant, we have that $C$ is isomorphic to an elliptic curve.  Thus by Siegel's theorem \cite[p. 353]{jhsdioph} $C$ has only finitely many integral points.
Therefore there are only finitely many $a_k$ with $d(a_k) = d(a_n)$.  Since there are only finitely many possible values for $d(a_n)$, we are done.

If $\Disc f = 0$, then $f = (x-k)^2$ for some $k \in \Z$.  We proceed as above, only this time we use that
$(a_{n-2}, c_n)$ gives an integral point on the curve defined by $f^2(x) = d(a_n)y^3,$ i.e.,
$$C : ((x-k)^2 - k)^2 = d(a_n)y^3.$$
Note that the curve $C' : (x-k)^2 - k = d(a_n)y^3$ has the same function field (over $\overline{\Q}$) as $C$, and is nonsingular.  Thus the genus of $C$ is the same as the genus of $C'$.  Since $k \neq 0$, we have that the discriminant of $(x-k)^2 - k$ is nonzero, and we use the same argument as above to show $C'$ has genus 1.  We then apply Siegel's theorem, also as above.  \done

{\bf Remark}: It would be a straightforward, if somewhat messy, task to extend the above argument to polynomials of higher degree.

\bibliographystyle{plain}

\end{document}